# Numerical Algorithm Development for Optimizing the Engine Stroke of Linear Generators

*Tulus, Syahril and Ahmad Kamal Ariffin*
*Department of Mechanical and Materials Engineering,*
*Universiti Kebangsaan Malaysia*
*Bangi, 43600 Selangor DE, Malaysia*

**Abstract**
This paper presents the results of a numerical algorithm development to optimize the engine strokes in a linear engine incorporating combustion and kickback. Due to the free piston reciprocal movement occurring in linear engines, the stroke of linear engines cannot be determined by using of fixed position. The reciprocating motion is namely the result of pressure acting on the piston and kickback. Kickback bore size is the main parameter that can influence the velocity of the motion. A numerical algorithm of a variational principle has been developed to optimize the kickback bore size. The resulting function of velocity is important to control the compression ratio and improve of engine thermal efficiency.

**Introduction**

Linear internal combustion engines may find application in the generation of electrical power without the need linear to rotary motion. The elimination of the connecting rod and crankshaft can improve the efficiency of the engine significantly and reduced weight and cost are added advantages. The operation of this engine is distinct from that of a conventional slider-crank mechanism engine, insofar as the motion of the two horizontally opposed pistons is not externally constrained [1]. This technology is advantageous because it is mechanically simpler and allows for a great deal more freedom in defining a piston motion profile, enabling the use of novel combustion regimes [2].

Many researchers have been designing various linear engines. Most of them have been developing linear engine using combustion at each cylinder end. Sandia National Laboratories has been investigating a new, integrated approach to generating electricity with ultra low emissions and very high efficiency for low power applications such as hybrid vehicles and portable generators. Their approach utilizes a free piston in a double-ended cylinder. Combustion occurs alternately at each cylinder end, with intake/exhaust processes accomplished through a two-stroke cycle [3].

The mean piston speed is and indicator of how well an engine handles load such as friction, inertia, and gas flow resistance. The mean piston speed is essentially a physical limitation and the operating speed depends upon the stroke. In a free piston configuration, the relationship between engine speed, operating conditions, and design parameters is much more complicated [4].

In the case that the combustion occurs only in one side of the linear engine, the problem is to manage the velocity of piston such that the maximum displacement does not exceed the maximum position. This paper is to present the results of numerical experiments using variational principle to get the optimal velocity and maximum position of piston.

**Problem design**

Consider a two stroke linear engine model as depicted in Fig. 1. To derive the governing equation, a priori, the following independent variables are specified. The left and right cylinder bore are $b_l$ and $b_r$ respectively. The maximum theoretical half-stroke length is $x_m$. The piston on a left to right stroke traverses from $-x_s$ to $+x_s$. The friction force required to move the slider is $F_f$. The quantity of heat added to one cylinder during a cycle is $Q_{in}$.

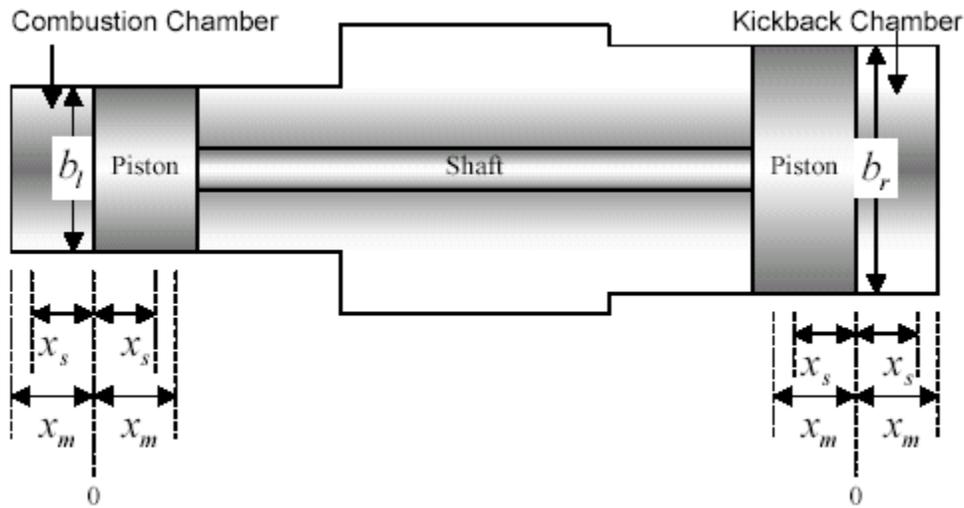

*Fig. 1. Simplified View of the Linear Engine*

Now let us consider a force balance for a left-to-right stroke, with $x$ positive in the left-to-right direction. The force balance can be written as

$$P_l(x)A_{Bl} - P_r(x)A_{Br} - F_f(x) = m_s \frac{d^2x}{dt^2} \quad (1)$$

where $P_l$ and $P_r$ are the instantaneous pressure in the left and right cylinder, respectively. By using an assumption of adiabatic compression and expansion behavior, the equation (Eq. 1) can be derived to get the following equation (Eq. 2).

$$P_{ml}\left(\frac{x_m}{x_m+x}\right)^n \frac{\pi b_l^2}{4} - P_{1r}\left(\frac{x_m+x_s}{x_m-x}\right)^n \frac{\pi b_r^2}{4} + Q_{in}(n-1)\frac{(x_m-x_s)^{n-1}}{(x_m+x)^n} - F_f(x) = m_s \frac{d^2x}{dt^2} \quad (2)$$

where $Q_{in}$ is a quantity of heat added during one cycle. This is a second order ordinary differential equation.

In case that the kickback bore will be determined, assume that the bore of kickback is scaled by the bore of the combustion cylinder. By letting $b_r = \lambda b_l$, the Eq. 2 can be written as follows

$$P_{ml}\left(\frac{x_m}{x_m+x}\right)^n \frac{\pi b_l^2}{4} - \lambda^2 P_{1r}\left(\frac{x_m+x_s}{x_m-x}\right)^n \frac{\pi b_l^2}{4} + Q_{in}(n-1)\frac{(x_m-x_s)^{n-1}}{(x_m+x)^n} - F_f(x) = m_s \frac{d^2 x}{dt^2} \quad (3)$$

The problem now is to solve the Eq. 3 in (0,T], with the initial conditions $x(0) = 0$, $\frac{dx}{dt}(0) = 0$, and $x_{max} = x_s$.

It is assumed that the bore scale $\lambda$ as a parameter and satisfies the condition

$$\lambda_{max} \geq \lambda \geq \lambda_{min} > 0 \quad (4)$$

where $\lambda_{max}$ and $\lambda_{min}$ are given positive constants.

For each value of $\lambda$ in Eq. 4, we can obtain the numerical solution of the Eq. 3 using numerical integration. But, in this problem, the value of $\lambda$ is unknown. Therefore variational methods are suitable to control the values of $\lambda$.

**Variational method and numerical algorithm**

The unknown parameter $\lambda$ is determine by minimizing the functional $J$ defined by

$$J(\lambda) := |x_{max}(\lambda) - x_s| \quad (5)$$

where $x_{max}$ is the maximum value of solution function $x$ as the velocity, $v = dx/dt$ is equal to zero. To minimize the functional $J$, the following procedure is used:

$\lambda_0$ is prescribed to initialize the bore scale. For $j = 0, 1, 2, \ldots$,

$$\lambda_{j+1} = \lambda_j + p_j s_j. \quad (6)$$

where $p_j$ is a search direction and $s_j$ is a search step. According to the definition (Eq. 5), it is possible that among the iteration, $\lambda_j$ either less or more than the appropriate $\lambda$. The search direction $p_j$ is obtained by

$$p_j = \begin{cases} 1 & \text{if } x_s \leq \max(\lambda_j) \\ -1 & \text{otherwise} \end{cases}$$

If $(x_s - \max(\lambda_j))(x_s - \max(\lambda_{j+1})) < 0$, then the search step is assigned as

$$s_{j+1} = s_j \frac{|x_s - \max(\lambda_{j+1})|}{|\max(\lambda_j) - \max(\lambda_{j+1})|}, \quad (7)$$

and otherwise

$$s_{j+1} = s_j. \quad (8)$$

The factor in the right hand side of the Eq. 7,

$$\frac{|x_s - \max(\lambda_{j+1})|}{|\max(\lambda_j) - \max(\lambda_{j+1})|} < 1.$$

Therefore, by using the Eqs. 7 and 8, the sequence of the values of $\lambda_{j+1}$ will tend to a limit.

**Results and Discussion**

In this section, a numerical experiment about our algorithm has been done using the following values of parameters: $P1_l$ = 225000 Pa, $P1_r$ = 120000 Pa, $Q_{in}$ = 18 Joule, $x_s$ = 0.0225 m, $b_l$ = 0.05 m, $F_f$ = 0, $n$ = 1.33. The maximum of positions for each value of $\lambda$ have been computed using numerical integration of the Eq. 3. By using the procedure above the initial bore scales of 1 and 2 have been chosen by assuming that the kickback bore must be larger than the combustion one. The iteration results are figured in Figs. 2 (a) and (b), respectively. The bore scale tends to 1.461, as the max $x$ tends to $x_s$. It means that the optimal piston velocity is achieved as we used a cylinder on the kickback side with bore size of 0.073 m. The function of optimal velocity and position are shown in Fig. 3.

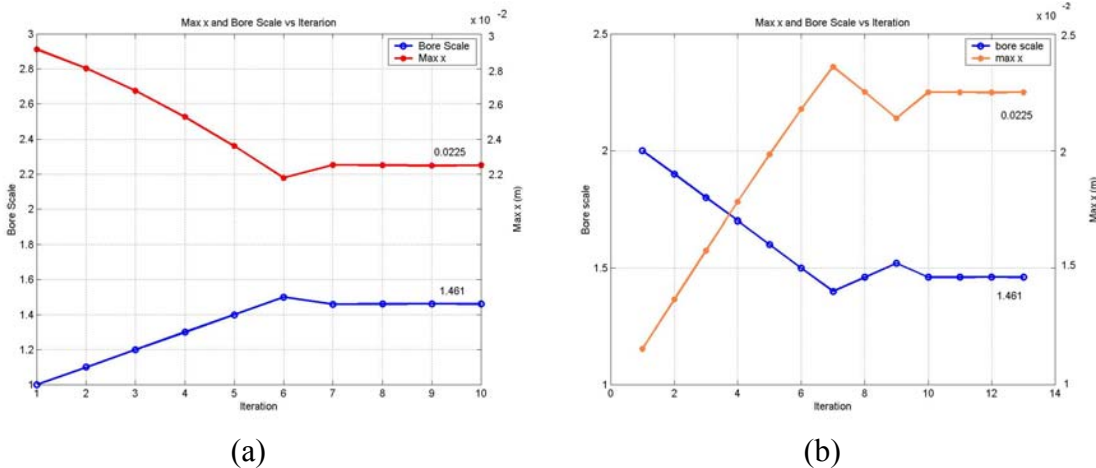

(a) (b)
***Fig. 2*** *Graphs of Max x and Bore Scale vs Iteration: (a) Initial bore scale is 1, (b) Initial bore scale is 2.*

**Concluding Remarks**

A numerical algorithm for the problem of bore scale identification of the cycle of linear two-stroke internal combustion engine has been presented. The direct variational method has been used to identify the unknown bore scale function. The numerical algorithm is based on minimizing the functional. The method used here converged fast after its values in one step across the required maximum displacement.

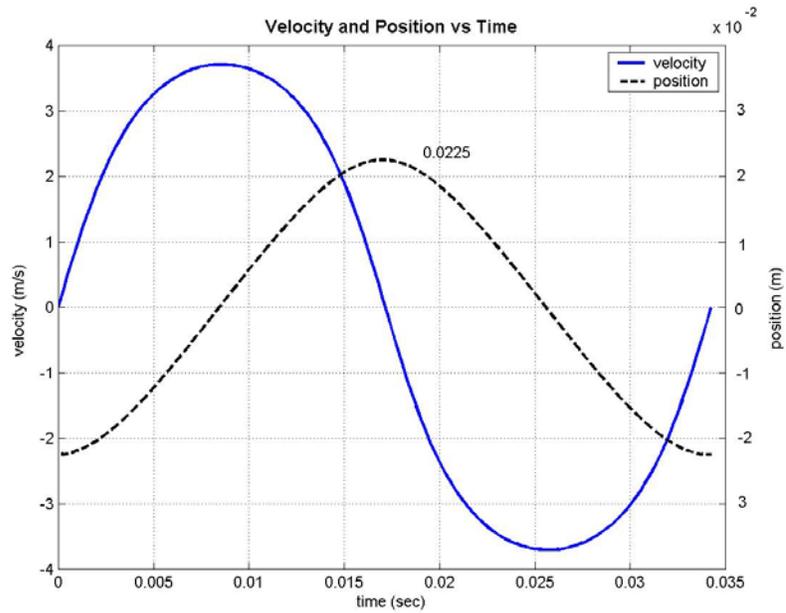

*Fig. 3*. Graph of the optimal velocity and position